\documentstyle[12pt,epsf]{article}
\def\fnote#1#2{\begingroup\def\thefootnote{#1}\footnote{#2}\addtocounter
{footnote}{-1}\endgroup}
\def\inbar{\vrule height1.5ex width.4pt depth0pt}
\def\IB{\relax{\rm I\kern-.18em B}}
\def\IC{\relax\,\hbox{$\inbar\kern-.3em{\rm C}$}}
\def\ID{\relax{\rm I\kern-.18em D}}
\def\IE{\relax{\rm I\kern-.18em E}}
\def\IF{\relax{\rm I\kern-.18em F}}
\def\IG{\relax\,\hbox{$\inbar\kern-.3em{\rm G}$}}
\def\IH{\relax{\rm I\kern-.18em H}}
\def\II{\relax{\rm I\kern-.18em I}}
\def\IK{\relax{\rm I\kern-.18em K}}
\def\IL{\relax{\rm I\kern-.18em L}}
\def\IM{\relax{\rm I\kern-.18em M}}
\def\IN{\relax{\rm I\kern-.18em N}}
\def\IO{\relax\,\hbox{$\inbar\kern-.3em{\rm O}$}}
\def\IP{\relax{\rm I\kern-.18em P}}
\def\IQ{\relax\,\hbox{$\inbar\kern-.3em{\rm Q}$}}
\def\IR{\relax{\rm I\kern-.18em R}}
\def\IT{\relax{\rm I\kern-.18em T}}
\def\ZZ{\relax{\sf Z\kern-.4em Z}}
\def\a{\alpha}   \def\b{\beta}     
   \def\k{\kappa}  
     \def\si{\sigma}

   \def\cD{{\cal D}}
\def\cF{{\cal F}}   
   
 \def\cO{{\cal O}}

 \def\bq{{\bar q}}

\def\fnote#1#2{\begingroup\def\thefootnote{#1}\footnote{#2}\addtocounter
{footnote}{-1}\endgroup}
\def\beq{\begin{equation}}  \def\eeq{\end{equation}}
\def\bea{\begin{eqnarray}}  \def\eea{\end{eqnarray}}
 \def\lleq#1{\label{#1}\eeq}

\def\tabroom{\hbox to0pt{\phantom{\Huge A}\hss}}
\def\notin{\ \hbox{{$\in$}\kern-.51em\hbox{/}}}

\def\del{\partial} 

\begin{document}
\hfill{alg-geom/9612012}
\vskip .01truein
\hfill{BONN--TH--96--13}
\vskip 1.4truein
\noindent
\centerline{{\bf SCALING BEHAVIOR ON THE SPACE OF CALABI--YAU MANIFOLDS}}

\vskip .8truein
\centerline{\sc R.Schimmrigk
                \fnote{\diamond}{netah@avzw02.physik.uni-bonn.de}
           }

\vskip .2truein
\centerline{\it Physikalisches Institut, Universit\"at Bonn}
\vskip .05truein
\centerline{\it Nussallee 12, 53115 Bonn}

\vskip 1.4truein
\centerline{\bf ABSTRACT}
\vskip .2truein

\noindent
Recent work is reviewed which suggests that certain universal quantities, 
defined for all Calabi--Yau manifolds, exhibit a specific behavior 
which is not present for general K\"ahler manifolds. 
The variables in question,  natural from a mathematical perspective,  
are of physical importance because they determine aspects of the 
low--energy physics in four dimensions, such as Yukawa couplings 
and threshold corrections.
It is shown that these quantities, evaluated on the complete class of 
Calabi--Yau hypersurfaces in weighted projective 4--space, exhibit scaling 
behavior with respect to a new scaling variable. 

\vskip .8truein

\centerline{\sc To appear in Mirror Symmetry II}

\renewcommand\thepage{}
\newpage

\baselineskip=18pt
\parskip=.2truein
\parindent=20pt
\pagenumbering{arabic}

\vskip .3truein
\noindent
{\bf Introduction.} 
A longstanding question in string theory concerns the structure of the space 
of groundstates. Much insight has been gained over the past
years into some of its salient properties for (2,2)--supersymmetric vacua 
\fnote{1}{In this paper the focus will be on the space of Calabi--Yau manifolds.}.
Most of these results, such as the discovery of mirror symmetry \cite{ms} 
and the fact that different Calabi--Yau manifolds are connected through 
singular configurations \cite{cdls88}, focus on the detailed geometry
of the moduli space of individual internal manifolds  
\fnote{2}{Or their field theoretic counterparts in the framework of
          Landau--Ginzburg theories or conformal field theories.}.
Relevant in this context is the generation--antigeneration structure of the
heterotic string, encoded in the two nontrivial Hodge numbers
$(h^{(1,1)},h^{(2,1)})$ of the Calabi--Yau space. These
topological numbers provide a first parametrization of the space of
string vacua and mirror symmetry, an operation that exchanges
these two numbers, provides a profound technique via which it is
possible to extract physical information about individual ground states.

What mirror symmetry has not provided, however, is a framework which leads 
to insight into the more detailed structure of the space of string vacua 
and which will enable us, eventually, to address problems like the vacuum 
degeneracy.  In the absence of a physical `theory of moduli space'
what is needed is some concrete information about this space, which we
might hope to be encoded in some numerical characteristics, much like 
data from an experiment.
Such numerical information is expected to be useful in providing a guideline 
as to what exactly it is that a future theory of moduli space has to explain.
In ref. \cite{rs94} it has been shown that such characteristics exist  
in the form of scaling exponents, thereby taking the first steps toward
a quantitative understanding of the moduli space of the supersymmetric
closed string. The scaling behavior defined by these exponents is 
specific to Calabi--Yau manifolds and does not exist for general projective 
varieties.

From a mathematical perspective the results of \cite{rs94} form part  
of a long tradition in algebraic geometry, which runs under the heading 
`geography of manifolds' and which reaches back at least as far as the last 
century.  As will become clear in the course of this review the analysis  
of \cite{rs94} is reminiscent of Noether's investigation of algebraic surfaces 
and represents a continuation of the work of Wilson \cite{pmhw94} in a slightly 
different context. 

Finally, we will have opportunity to comment on recent results \cite{ko91} 
regarding the classification of Calabi--Yau manifolds using the Fujita 
$\Delta$--index. 
 
\vskip .1truein
\noindent
{\bf The Variables.} 
The quantities which will be focused on in the following are based on
the existence of a universal structure which exists on every Calabi--Yau
manifold by virtue of the fact that they are projective, {\it i.e.} they are
described by equations in a particular type of compact, complex, K\"ahler
manifold, so--called projective spaces.  On such spaces there exists a natural
structure, the so--called hyperplane bundle, denoted by $L$, which can
be restricted to the Calabi--Yau manifold embedded in the ambient projective
space. This bundle is of interest for physics because its first Chern class
$c_1(L)$ leads to a universal massless mode, present in every
Calabi--Yau manifold, which parametrizes one of the antigenerations
observed in four dimensions. Moreover, on a three--complex dimensional
Calabi--Yau manifold, the degree of this line bundle, defined as
\beq
L^3 \equiv \int_M c_1^3(L),
\eeq
affords an interpretation as a Yukawa coupling $\k_L \equiv L^3$ of the
antigeneration associated to $L$. Thus $L^3$ determines
the strength of the Yukawa coupling of the corresponding
four--dimensional mode. A second characteristic of any line bundle is the
number of its global functions, denoted by $h_L$. For every Calabi--Yau manifold
the line bundle $L$ thus leads to a pair of quantities $(L^3,h_L)$
which can be viewed as coordinates on the space of Calabi--Yau manifolds.
The question arises whether these new coordinates allow us to draw
some type of `phase diagram' for the spaces of interest, {\it i.e.}
whether they provide a new sort of structure theorem for string vacua, 
furnishing a new `geographical map'.

A final quantity of interest which is induced by the line bundle $L$ is
\beq
L\cdot c_2 \equiv \int_M c_1(L)\wedge c_2(M),
\lleq{c2}
where $c_2(M)$ is the second Chern class of the Calabi--Yau 3--fold.
The physical interpretation of this quantity has been uncovered
in ref. \cite{bcov93}, where it was shown that the generalized index
\beq
\cF =\frac{1}{2} \int \frac{d^2 \tau}{{\rm Im}~\tau}~
         Tr\left[(-1)^F F_L F_R q^{H_L} \bq^{H_R}\right],
\eeq
introduced in \cite{cfiv92}, contains information about the threshold
corrections of the gauge couplings in string theory (see also \cite{agn92}) 
and that it
is the key for the understanding of quantum mirror symmetry at one loop.
Here the integral is over the fundamental domain of the moduli space of the
torus, $F_{L,R}$ denote the left and right fermion numbers and the trace is
over the Ramond sector for both the left-- and right--movers
\fnote{3}{The contribution of the ground states of the supersymmetric Ramond
     sector to $\cF$  has to be deleted in order for the integral to
converge.}.
In lowest order this index essentially reduces to
$\frac{1}{24} \int_{M} K\wedge c_2(M)$, where $K$ is the K\"ahler form of
the manifold. Hence the second Chern class evaluated on $L$, (\ref{c2}),
defines the universal contribution $\cF^{\uparrow}_L \equiv L\cdot c_2/24$
to the large volume limit of this partition function.

A useful fact is that the three variables $L^3,h_L,L\cdot c_2$ are not all 
independent because of the theorem of Hirzebruch--Riemann--Roch, which relates
the Euler number defined by the counting problem   
\beq
\chi(M,L) = \sum_{p=0}^n (-1)^p {\rm dim~H}^p(M,L) 
\eeq
to an integral over the manifold
\beq
\chi(M,L) = \int_M e^{c_1(L)} \wedge {\rm Td}(M),
\eeq
where Td$(M)$ is the Todd class of the manifold. For Calabi--Yau 
3--folds this leads to 
$\chi(M,L) = \frac{1}{6}L^3+\frac{1}{12}L\cdot c_2$  and hence 
\beq
L^3 = 6h_L -\frac{1}{2} L\cdot c_2.
\eeq

Finally, we quote, for future use, the definition of the Fujita index 
\cite{tf90}  which can be defined for a variety of arbitrary dimension $n$ 
as 
\beq
\Delta = n + L^n - h_L. 
\eeq
This index has proven useful in the classification of polarized varieties 
in particular when it assumes small values.

\vskip .1truein
\noindent
{\bf The Class.} 
A class of Calabi--Yau spaces that is particularly amenable to an analysis
in terms of the variables $L^3$, $L\cdot c_2$, and $h_L$ is
furnished by hypersurfaces embedded in weighted projective space
$\IP_{(k_1,...,k_5)}$. The complete class of such manifolds, consisting of
7,555 configurations, has been
constructed in \cite{ks94,krsk92}. The natural candidate for
a line bundle on such spaces is the pullback of the weighted form
of the hyperplane bundle  $L \equiv \cO^{(k)}_{\IP_{(k_1,...,k_5)}}$
defined on the weighted ambient space, where  \cite{d93,bk93}
\beq
k=lcm\{~\{gcd(k_i,k_j)|~i,j=1,..,5; i\neq j\} \cup
\{k_i| k_i ~ {\rm does~not~divide}~ \sum_{i=1}^5 k_i \} \}.
\eeq
The restriction of $L$ from the weighted ambient space to the embedded
Calabi--Yau manifold $M$ induces an antigeneration which will also be
denoted by $L$. The Yukawa coupling $\k_L$ of this antigeneration leads,
in the large volume limit, to the expression of the degree of $L$ 
\beq
L^3 = \int_M (j^*(c_1( \cO^{(k)}_{\IP_{(k_1,...,k_5)}}))^3
      =\left(\frac{\sum_{i=1}^5 k_i}{\prod_{i=1}^5 k_i}\right)~k^3.
\eeq
The number of global functions of the bundles
$L=\cO^{(k)}_{\IP_{(k_1,...,k_5)}}$ can be obtained as
\beq
h_L = \frac{1}{k!} \frac{\del^k}{\del t^k}
 \left(\frac{\left(1-t^{\sum_{i=1}^5 k_i}\right)}{\prod_{i=1}^5
   \left(1-t^{k_i}\right)} \right)~\rule[-4mm]{.1mm}{10mm}_{~t=0} 
\eeq
and therefore we can determine the Fujita index and also, via the 
Hirzebruch--Riemann--Roch--relation 
the second Chern class evaluated on the hyperplane bundle.
In the following these variables will be used to explore the structure
of the moduli space of Calabi--Yau manifolds.  

\noindent 
{\bf The Results.}

{\it The Yukawa coupling:} For arbitrary projective manifolds the two numbers
$L^3$ and $h_L$ are independent, completely uncorrelated, quantities.
Hence in a diagram of these numbers we might expect a random distribution of
data points, with no obvious pattern. The actual computation for 
Calabi--Yau manifolds, the results
of which are shown in Figure 1, however uncovers an unexpected simplicity.

\vskip .3truein
\centerline{ \epsfbox{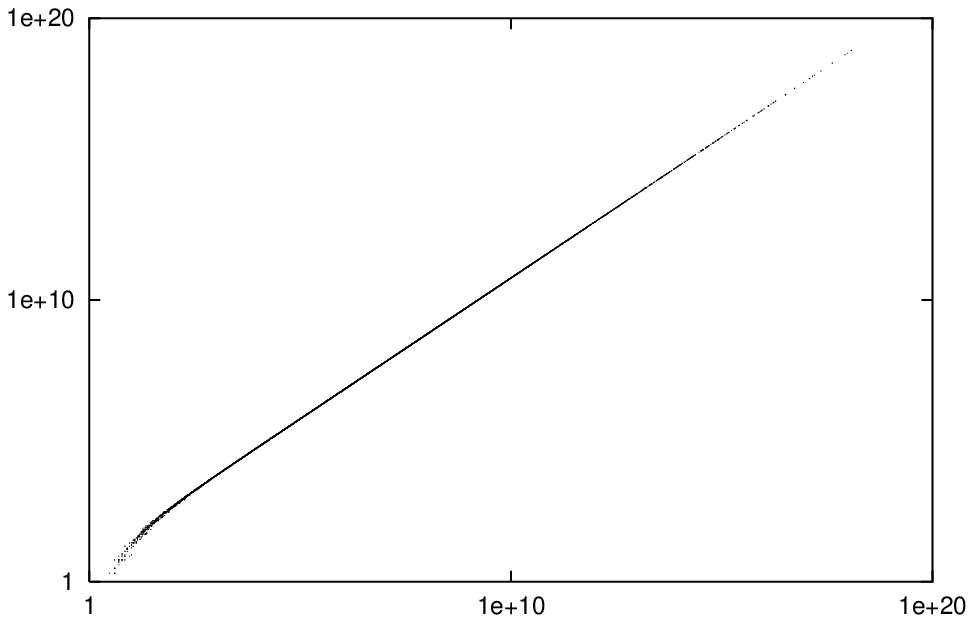}}

\vskip .2truein
\centerline{{\bf Figure 1:} {\it The degree $L^3$ versus $h_L$.}}

The first thing to notice is the emergence of a well defined boundary
which is described by a simple relation: all varieties in the class of
Calabi--Yau 3--folds embedded in weighted $\IP_{(k_1,...,k_5)}$
satisfy the inequality
\fnote{4}{The weaker limit $\k_L < 6h_L$ follows from the relation of
          Hirzebruch--Riemann--Roch and the positivity of
          the second Chern class.}
\beq
L^3 \leq  2 \left(3h_L-8\right).
\lleq{yukasyl}
This is in contrast to mirror symmetry where the boundary of the
mirror diagram (see the first reference in \cite{ms}) does not yield to
such an easy description, even though it is also well defined.

What is rather unexpected is that as the dimensions $h_L$ increase for
the individual Calabi--Yau hypersurfaces the Yukawa couplings approach the
upper limit
$L^3_{as} \equiv 2(3h_L-8)$  
very quickly, with very little scattering. Furthermore they do so according
to a power law.  This can be seen most clearly from Figure 2, which contains
the distances $\cD \equiv (L^3_{as}/L^3-1)$ 
of the degree of $L$  from $L^3_{as}$ as a function of $h_L$.

\vskip .3truein

\centerline{ \epsfbox{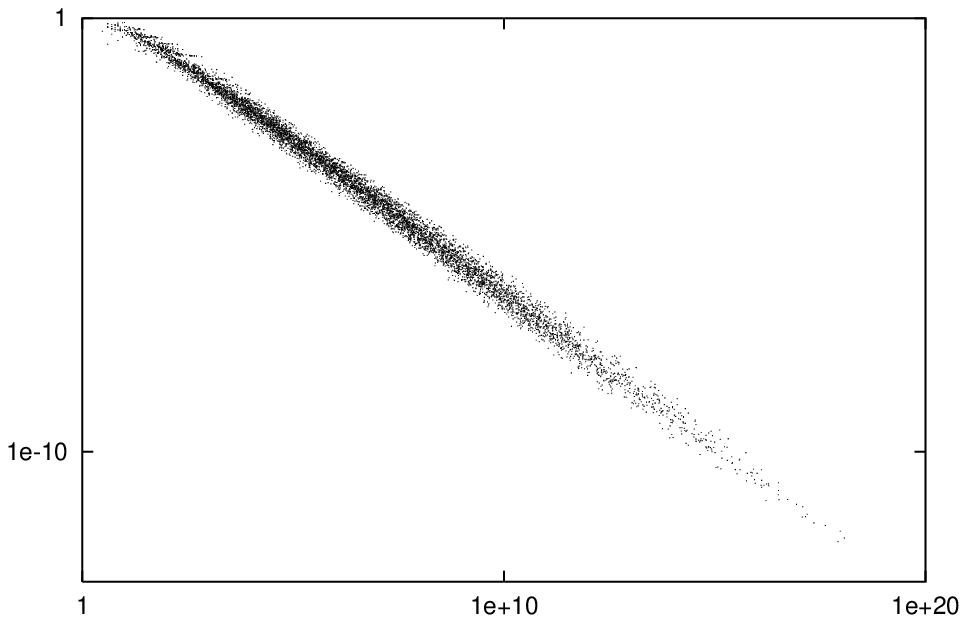}}

\vskip .1truein
\centerline{{\bf Figure 2:} {\it $\cD$ versus $h_L$.}}

The result shows clearly that the distance variable $\cD$ exhibits a scaling behavior
with approximate exponents $ A=5$  and $\a = 0.7$, 
determined by a least square fit. Solving this scaling relation leads to the 
functional dependence of the degree of the hyperplane bundle on $h_L$.  
These  results can be summarized more precisely in the following 

\noindent
{\bf Theorem 1.} {\it All varieties in the class of CY 3--folds embedded in 
   weighted projective 4-space $\IP_{(k_1,...,k_5)}$ satisfy the inequality }
\beq
 L^3 \leq 2 \left(3h_L-8\right).
\eeq
{\it The distance $\cD \equiv \frac{2(3h_L-8)}{L^3} - 1$
from the asymptotic line $L^3_{as} \equiv 2(3h_L-8)$
scales like} 
 \beq
  D \sim A (h_L)^{-\a}
  \lleq{discale}
{\it with $A\sim 5$ and $\a \sim 0.7$. Thus the functional dependence of the 
 degree on $h_L$ is given by }
\beq
L^3  \sim  \frac{L^3_{as}}{1+A (h_L)^{-\a}}   
\lleq{yukpow}

\noindent
{\bf Remarks.}
\begin{itemize}
\item Physically therefore the Yukawa coupling $\k_L\equiv L^3$ behaves like an 
      order parameter with the dimension $h_L$ assuming the r\^{o}le of the scaling 
     variable.
\item The degree $L^3$  and the dimension $h_L$  take values over an enormously wide 
       range of values, unprecedented in Calabi--Yau physics, spanning 18 orders of 
      magnitude. This will propagate via the Hirzebruch--Riemann--Roch relation to 
     $L\cdot c_2$ and lead to an enormous variation for the threshold corrections 
    of these string vacua (see remark below). 
\item  The power law behavior (\ref{yukpow}) for the coupling provides some measure 
    of simplicity of Calabi--Yau manifolds: the larger the value of the
scaling variable $h_L$, the more accurate the prediction of the Yukawa
coupling via (\ref{yukpow}). Since, very roughly, the magnitude
of the scaling variable is determined by the degree of complexity of the
singular sets of the varieties that need to be resolved, spaces with
more complicated resolution geometry are in fact simpler when it comes
to scaling. This phenomenon is rather different from the conventional
point of view which considers smooth projective hypersurfaces as the
simplest type of manifolds.  
\item Finally, it is easy to see that the scaling
      behavior (\ref{discale}), (\ref{yukpow}) is a characteristic of Calabi--Yau
manifolds, not of general projective 3--folds, not to mention general
3--folds: consider the hypersurfaces $\IP_4[d]$ of degree $d$
embedded in ordinary projective 4--space. For this (infinite) class of
manifolds one finds $L^3 =d, h_L=5, L \cdot c_2 = d(d^2-5d+10).$
\end{itemize}

\vskip .1truein
{\it The second Chern class.} 
Over the last years a confluence of ideas has occured in physics and mathematics 
toward a further physically important quantity. Wilson \cite{pmhw94},
on the one hand, has reemphasized the importance of the linear form on the
Picard group, defined by the second Chern class $c_2$, for purposes of
classification, a fact first recognized by Wall.
Bershadsky, Cecotti, Ooguri, and Vafa, on the other, have found that the
generalized index of \cite{cfiv92} reduces in the large volume limit
essentially to $\cF_L^{\uparrow}=L\cdot c_2/24$, and have shown that this
quantity arises in the 1--loop quantum mirror expansion
of the Ray--Singer torsion, as well as in the threshold corrections of the 
gauge couplings (see also \cite{agn92}). Thus the pair
$(L^3, L\cdot c_2)$ provides two coordinates on the moduli space
of Calabi--Yau manifolds with immediate physical meaning.
To investigate a possible relation between these quantities
is of particular interest since they arise in two different perturbation
expansions of the string -- the $\si$--model expansion at string tree level
and the string loop expansion.

Wilson \cite{pmhw94} has observed that for the degree  $L^3$ and
$L\cdot c_2$ an inequality $L\cdot c_2 \leq 10 L^3$ is obtained
for models whose Fujita index $\Delta$ is greater than two. The analysis
of the present class of weighted Calabi--Yau hypersurfaces shows that for
all but 11 spaces the Fujita index is always larger than two. 
Hence Wilson's observation shows that except for those sixteen spaces the
inequality above does hold in this class, i.e. in physics notation 
\beq
\frac{\cF^{\uparrow}_L}{\k_L} \leq \frac{5}{12}.
\lleq{c2uplim}
No further information about the population below this upper limit is known.

The unexpected consequence of the analysis above of all Calabi--Yau
hypersurfaces in weighted projective 4--space
is that the vast majority of these vacua does not come even
close to the upper limit (\ref{c2uplim}) but instead takes values in a narrow
region far below that line.  For large values of $L\cdot $ the
Yukawa couplings lie more than ten orders of magnitude below the upper
limit provided by Wilson's observation.  Furthermore the results, as shown
in Figure 3, uncover a definite correlation between the degree of $L$ and
$L\cdot c_2$.

\vskip .3truein
\centerline{ \epsfbox {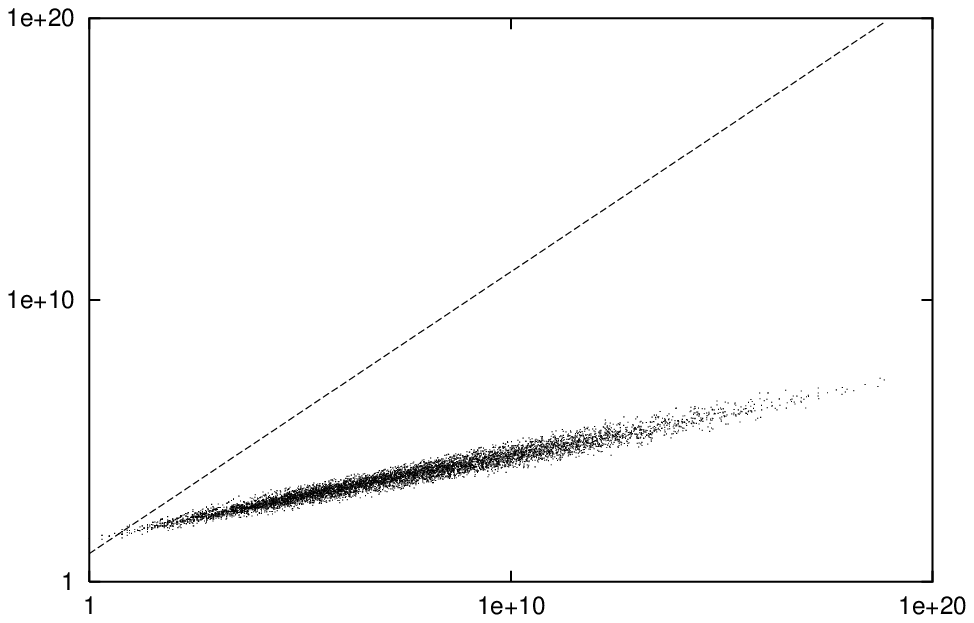}}

\vskip .2truein

\noindent
\centerline{{\bf Figure 3:}{\it~~ $L\cdot c_2$ versus $L^3$.
              The dotted line represents Wilson's upper limit.}}

An analysis similar to the one discussed previously leads to the following 

\noindent
{\bf Theorem 2.} {\it In the class of Calabi--Yau hypersurfaces in weighted 
projective space $\IP_{(k_1,...,k_5)}$ the product $L\cdot c_2$ scales like }  
\beq
L \cdot c_2 \sim B ~\k_L^\beta
\eeq
{\it with approximate values }
\beq
B =36,~~~~\beta=0.29.
\eeq

\noindent
{\bf Remark.} As mentioned above, the range of $L \cdot c_2$ is huge, resulting in 
      vast range of low--energy effective GUT mass scales
     \cite{bcov93}
      \beq
      M_{\rm GUT}^{\rm eff} =
      M_{\rm GUT}~ exp\left( \begin{small} \frac{1}{b}\int_M K\wedge
          c_2(M)\end{small}\right), 
     \eeq
    where $b=54+3(h^{(1,1)}+h^{(2,1)}$. 

\vskip .1truein
{\it The Fujita index:} 
From the definition of this index it is clear that it will behave precisely 
like the degree of the hyperplane bundle when viewed as a function of $h_L$. 

Of particular interest to mathematicians have been manifolds with low values 
of the $\Delta$--genus. In \cite{ko91} all smooth Calabi--Yau manifolds with 
$\Delta \leq 2$ have been classified. In the set of all Calabi--Yau hypersurfaces 
in weighted projective space there are only 4 smooth members. The results of 
\cite{rs94} lead to an immediate generalization of Oguiso's analysis and show 
that the class of manifolds with low Fujita index broadens considerably when 
resolutions are admitted. 
Table 1 contains the weighted configurations with $\Delta \leq 2$ together 
with some of their pertinent properties. 

\vskip .1truein
\begin{center}
\begin{small}
\begin{tabular}{|| l c c c c c c r ||}
\hline
\hline
Manifold &$\Delta$ &$L^3$ &$h_L$ &$L\cdot c_2$ &Euler number &$h^{(1,1)}$
&$h^{(2,1)}$\tabroom
\\
\hline
\hline
$\IP_{(1,1,1,2,5)}[10]$ &1 &1 &3  &34  &--288 &1 &145  \tabroom \\
$\IP_{(1,2,2,2,7)}[14]$ &1 &2 &4  &44  &--240 &2 &122  \tabroom \\
$\IP_{(1,1,1,1,4)}[8]$  &1 &2 &4  &44  &--296 &1 &149  \tabroom \\
$\IP_{(1,2,2,3,4)}[12]$ &2 &2 &3  &32  &--144 &2 &74   \tabroom \\
$\IP_{(1,3,3,3,5)}[15]$ &2 &3 &4  &42  &--144 &3 &75   \tabroom \\
$\IP_{(1,2,3,3,9)}[18]$ &2 &3 &4  &42  &--192 &3 &99   \tabroom \\
$\IP_{(1,1,1,1,2)}[6]$  &2 &3 &4  &42  &--204 &1 &103  \tabroom \\
$\IP_{(1,1,2,2,6)}[12]$ &2 &4 &5  &52  &--252 &2 &128  \tabroom \\
$\IP_{(1,2,3,12,18)}[36]$ &2 &6 &7 &72  &--360 &5 &185 \tabroom \\
$\IP_{(1,1,2,8,12)}[24]$  &2 &8 &9 &92  &--480 &3 &243 \tabroom \\
$\IP_{(1,1,1,6,9)}[18]$ &2 &9 &10 &102  &--540 &2 &272  \tabroom \\
\hline
\hline
\end{tabular}
\end{small}
\end{center}

\centerline{{\bf Table 1:} {\it All Calabi--Yau manifolds embedded in weighted $\IP_4$
                  with $\Delta \leq 2$. }}

This indicates that there should be still more to explore once one is able to 
probe more general weighted complete intersections of hypersurfaces in toric 
varieties and their generalizations.

\noindent
{\bf In summary}, the results presented above suggest that the degree  
$L^3$ and the universal lowest order part $L\cdot c_2$  of the free energy should
be interpreted as order parameters for which the number $h_L$ of global
functions behaves like a scaling variable. The scaling phenomena described
above provide specific power laws with a particular set of exponents that
any future physical theory of moduli space of the heterotic string will have
to explain. Thus our scaling laws furnish the first quantitative insight into
the structure of the string configuration space.

The present situation is not without historical precedent in physics: it is
reminiscent of Bjorken's results many years ago when he argued that a certain
scaling behavior should manifest itself in deep inelastic electron--nucleon
scattering. Our scaling can in fact be presented in a form which is
very similar to the original data of these historical results.
Define a function $(h_L)^{-1} \k_L$ and consider its dependence
on the scaling variable $h_L$. The result, which is displayed in Figure 4,
 looks very much like the experimental results of the SLAC group for
$\nu W_2(\nu,q^2)$.

\vskip .2truein
\centerline{ \epsfbox {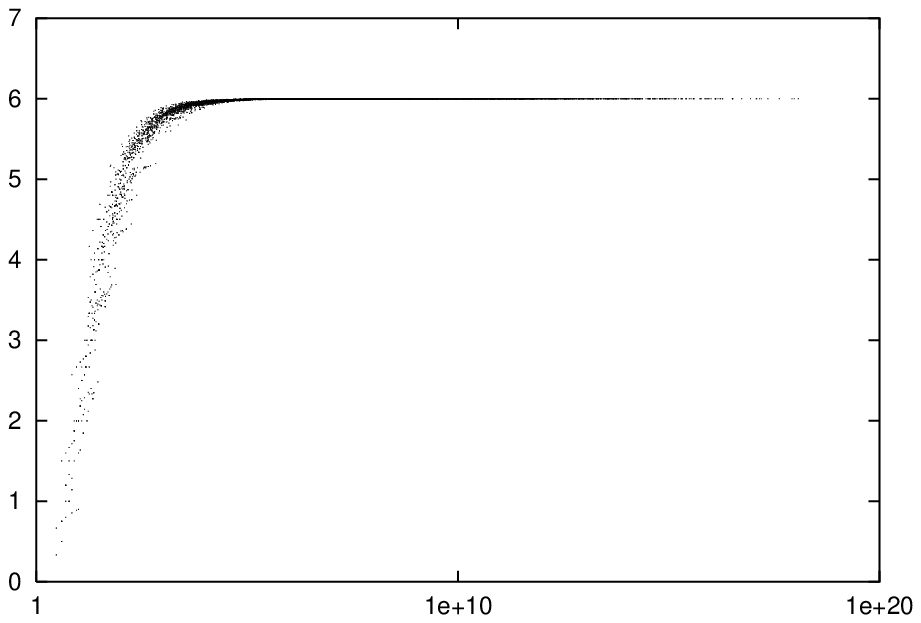}}

\vskip .2truein

\noindent
\centerline{{\bf Figure 4:}{\it~~ $(h_L)^{-1} \k_L$ versus $h_L$.}}

\noindent
What is needed, then, is the analog of Feynman's quark--parton picture, or
perhaps rather the analog of the QCD explanation of, say,
the momentum fraction of the proton carried by the quarks.

Tantalizing questions remain: An immediate one is
whether the exponents $\a,\b$ are universal or specific to our class of
Calabi--Yau hypersurfaces in weighted projective 4--space. There are two
natural avenues for further exploration. The simpler of the two is provided
by the toric framework of ref. \cite{b94} which has been shown \cite{cok94}
to include all Calabi--Yau hypersurfaces in weighted projective 4--space.
The second direction involves a natural generalization
of the framework of Calabi--Yau compactification, motivated by
mirror symmetry.  In \cite{rs93} a construction has been introduced which
provides an embedding of Calabi--Yau spaces into a larger class of a
special type of Fano manifolds, varieties whose  positive first Chern
class satisfies a particular quantization condition. It was shown in
\cite{rs93} that even though these special
Fano manifolds are not consistent ground states themselves they
do contain information about (2,2) supersymmetric string vacua.
It was furthermore demonstrated in ref. \cite{cod93} that such manifolds also
lead to the correct behavior of the Yukawa couplings and in ref. \cite{bb94}
that the framework of \cite{rs93} lends itself to a toric analysis,
generalizing the toric mirror construction of \cite{b94}. A natural question is
whether our scaling behavior persists in these more general settings as well.
Finally, the computations presented here are of lowest
order in perturbation theory and an interesting problem is to explore
the consequences of including higher order corrections.

\noindent
\section*{Acknowledgment}
It is a pleasure to thank Philip Candelas, Dimitrios Dais,
Xenia de la Ossa, Ed Derrick, Michael Flohr, Ariane Frey,
Jerry Hinnefeld, Vadim Kaplunovsky, Jack Morse, Werner Nahm, Steve Shore,
and especially Andreas Honecker, Monika Lynker and Katrin Wendland
for discussions.
I'm grateful to the Theory Group at the University of Texas at
Austin, the Department of Physics at Indiana University at South Bend,
and Simulated Realities Inc., Austin, TX for hospitality.

\end{document}